\documentclass[12pt]{article}
\usepackage{graphics, latexsym}
\usepackage{graphicx}
\usepackage{setspace}
\usepackage{amssymb}
\usepackage{amsmath}
\usepackage{float}
\usepackage{authblk}

\begin{document}
\title{Optimal Sequential Tests for Monitoring Changes in the Distribution of Finite Observation Sequences}
\author[1]{Dong Han}
\author[2]{Fugee Tsung}
\author[1]{Jinguo Xian}
\affil[1]{Department of Statistics, Shanghai Jiao Tong University, China.}
\affil[2]{Department of Industrial Engineering and Logistics Management,
Hong Kong University of Science and Technology, Hong Kong.}
\maketitle
\begin{abstract}
 This article develops a method to construct the optimal sequential test for monitoring the changes in the distribution of finite observation sequences with a general dependence structure. This method allows us to prove that different optimal sequential tests can be constructed for different performance measures of detection delay times.  We also provide a formula to calculate the value of the generalized out-of-control average run length for every optimal sequential test. Moreover, we show that there is an equivalent optimal control limit which does not depend on the test statistic directly when the post-change conditional densities (probabilities)  of the observation sequences  do not depend on the change time.
\end{abstract}

\renewcommand{\thefootnote}{\fnsymbol{footnote}}
\footnotetext{
Supported by RGC Competitive Earmarked Research Grants, National Basic Research Program of China
(973 Program, 2015CB856004) and National Natural Science Foundation of China (11531001)
\newline $\,\, MSC\, 2010\,subject\,classification $. Primary 62L10;
secondary 62L15}

\section{Introduction}

One of the basic problems in statistical process control (SPC) is designing an effective sequential test (or a control chart), as proposed by Shewhart (1931), to detect possible changes at some instant (change-point) in the behavior of a series of sequential observations. The objective is to raise an alarm as soon as a change occurs, while keeping the rate of false alarms to an acceptable level.  Detecting abrupt changes in a stochastic system quickly without exceeding a specified false alarm rate is an important issue not only in industrial quality and process control applications, but also in non-industrial processes (Bersimis  \emph{et al}. 2018), biology (Siegmund 2013), clinical trials and public-health (Woodall 2006,  Chen and Baron 2014, Rigdon and Fricker 2015) , econometrics
and financial surveillance (Fris\'{e}n 2009),  graph and network data ( Akoglu \emph{et al}. 2015,  Woodall \emph{et al}. 2017,  Hosseini and Noorossana 2018), etc.

A great variety of sequential tests have been proposed, developed and applied to detect changes in the distribution of sequential observations quickly in various fields; see, for example, Siegmund (1985),  Basseville and Nikiforov (1993), Lai (1995, 2001), Stoumbos \emph{et al}. (2000), Chakraborti \emph{et al}. (2001),  Bersimis \emph{et al}. (2007),  Montgomery (2009),  Poor and Hadjiliadis (2009), Woodall and Montgomery (2014),   Qiu (2014) and Tartakovsky \emph{et al}. (2015). This raises two questions: What is the optimal sequential test?   How do we design or construct an optimal sequential test?

First, we recall the main results of the known optimal sequential tests. A sequential test $T^*$ is called to be optimal for detecting changes in the distribution if the average value of some detection delay time $(T-k+1)^+$ of $T^*$ for  all possible change time $ k\geq 1$  is the smallest of all of the sequential tests $T$ with a given probability of false alarm that is no greater than  a preset level ( or with a given false alarm rate that is no less than a given value), where  $x^+=\max\{0, x\}$. In the literature, there are four main kinds of optimal sequential tests: the Shiryaev (1963, 1978, P.193-200) test $T_S(c_1)$, two SLR (sum of the log likelihood ratio) tests  $T_{SLR_1}(c_2)$ (Chow, Robbins and Siegmund 1971, P.108) and $T_{SLR_2}(c_3)$ (Fris\'{e}n 2003), the CUSUM test $T_C(c_4)$ (Page 1954, Moustakides 1986)  and the Shiryaev-Roberts test $T^r_{SR}(c_5)$ (Polunchenko and Tartakovsky 2010), where the five positive numbers $c_i>0,$ $1\leq i\leq 5$,  denote the five constant control limits or the threshold limits. It can be seen that to prove the optimality of the tests above we need the assumption that  there is an infinite independent observation sequences $\{X_n, \, n\geq 1\}$ or the corresponding test statistics are Markov sequences.

In fact, it is not realistic for us to have an infinite observation sequences, that is, people can only obtain finite observation sequences in reality.  For example, consider a production line that produces one product per minute. If the production line works eight hours a day, then  the number of products or observations per day is $N=480$. Our task is to design or construct an effect test for detecting  whether the 480 observations  (usually not independent) are abnormal in real-time. However, when we only have $N$ finite independent observation sequences $\{X_n, \, 1\leq  n\leq N\}$ ($N\geq 2$), all the five optimal sequential tests mentioned above will become no longer optimal, that is,  the five tests: $T_{S, N}(c_1)=\min\{T_S(c_1), N+1\}$, $T_{SLR_1, N}(c_2)=\min\{T_{SLR_1}(c_2), N+1\}$, $T_{SLR_2, N}(c_3)=\min\{T_{SLR_2}(c_3), N+1\}$, $T_{C, N}(c_4)=\min\{T_C(c_4), N+1\}$ and $T^r_{SR, N}(c_5)=\min\{T^r_{SR}(c_5), N+1\}$,  will be no longer optimal for finite independent observation sequences (see Corollary 3, Remark 5 in Section 3). Hence, how to construct an optimal sequential test for finite observation sequences will become very important.  However, since Shewhart (1931) proposed a control chart method for sequential testing, there has been little progress in constructing and proving the optimal sequential test for finite  observation sequences with a general dependence structure. The main purpose of this study is to try to solve this problem.

In this paper,  based on Chow-Robbins-Siegmund's work (1971, Chaper 3)  we develop a method to construct various optimal sequential tests under different performance measures of detection delay times for detecting the change in probability distribution of  finite observation sequences. Moreover,  we  find a formula to calculate the value of the generalized out-of-control average run length for each optimal test and obtain an equivalent optimal control limit which may not depend on the test statistic directly.  As a corollary of the above conclusion,  the five tests $T_{S, N}(c_1), T_{SLR_1, N}(c_2), T_{SLR_2, N}(c_3), T_{C, N}(c_4)$ and $T^r_{SR, N}(c_5)$ can be still optimal for detecting the change in finite observation sequences, if their constant control limits $c_i, 1\leq i\leq 5,$ are replaced by the corresponding so-called optimal dynamic control limits respectively (see Corollary 1 in Section 2.2).

The rest of this paper is organized as follows. Section 2.1 presents a generalized Shiryaev's measure to evaluate how well a sequential test performs to detect changes in the distribution of finite observation sequences. Section 2.2 constructs the  optimal sequential test and gives the formula for calculating the generalized out-of-control average run length. The equivalent optimal control limit is presented and proved in Section 3. The detection performance of two optimal tests is illustrated by comparison and analysis of the numerical simulations for 60 observations in Section 4. Section 5 provides some concluding remarks. Proofs of the theorems are given in the Appendix.

\section{Optimal sequential tests for finite observations}

In this section, we first present the performance measure and optimization criterion, then construct the optimal sequential tests.

Consider finite observations, $X_1, X_2, ..., X_N$.  Without loss of generality, we assume $N\geq 2$. Let $\tau=k$ ($ 1\leq k \leq N $) be the change-point.  Let  $p_{0}(x_0, x_{1},...,
x_N)$ and $p_{k}(x_0, x_{1},..., x_N)$ be the pre-change and post-change joint probability densities respectively. Denote the post-change joint probability
distribution and  the expectation by $\textbf{P}_{k}$ and  $\textbf{E}_{k}$ respectively for $1\leq k \leq N$.  When $\tau >N$, i.e., a change never occurs in $N$ observations $X_1, X_2, ..., X_N$, the probability
distribution and the expectation are denoted by $\textbf{P}_{0}$ and $\textbf{E}_{0}$ respectively for all observations $X_0, X_1, X_2, ...,X_N$ with the pre-change joint probability density $p_{0}(x_0, x_{1},...,x_N)$.
 Moreover, when the observation sequence $\{X_n, 0 \leq n\leq N\}$ takes discrete values, the above joint probability densities and the conditional probability densities will be considered as joint  probability distributions and the conditional probability distributions taking the discrete values.

In order to construct the optimal sequential tests in Section 2.2, we assume that the following likelihood ratio of the post-change conditional probability density to the pre-change conditional probability density, $\Lambda^{(k)}_j$, satisfies
\begin{eqnarray}
\Lambda^{(k)}_j=\frac{p^{(k)}_{\textbf{1}j}(X_j|X_{j-1}, ..., X_0)}{p_{\textbf{0}j}(X_j|X_{j-1}, ..., X_0)}<\infty \,\,\,( a.s. \textbf{P}_{0} )
\end{eqnarray}
and has no atoms with respect to $\textbf{P}_{0}$ for $1\leq k\leq N$ and $k\leq j\leq N$, where $p_{\textbf{0}j}(x_j|x_{j-1}, ..., x_0)$ for $1\leq j \leq N$ and $p^{(k)}_{\textbf{1}j}(x_j|x_{j-1}, ..., x_0)$  for $1\leq k \leq N$, $k\leq j\leq N,$ denote the pre-change and post-change conditional probability densities, respectively,  and the notation $(k)$ in $p^{(k)}_{\textbf{1}j}$ denotes that the post-change conditional probability densities $ p^{(k)}_{\textbf{1}j}$ rely on the change-point $k$  for $k\leq j\leq N$.  If $\Lambda^{(k)}_j=\Lambda_j$ for $1\leq k\leq i\leq N$, it means that the post-change conditional densities (probabilities)  of the observation sequences do not depend on the change-point.

\subsection{ Performance measures of sequential tests}
Let $T \in \mathfrak{T}_N$ be a sequential test, where $\mathfrak{T}_N$ is a set of all the sequential tests satisfying $1\leq T\leq N+1$ and $ \{T\leq n\} \in  \mathfrak{F}_{n}=\sigma\{X_{j}, 0\leq j\leq n\}$ for $1\leq n\leq N$. Let $W=\{w_j, 1\leq j\leq N+1\}$ and $V=\{v_j, 1\leq j\leq N+1\}$ be two series of nonnegative random variables satisfying $w_k, v_k \in \mathfrak{F}_{k-1}$ for $1\leq k\leq N+1$. Denote the indicator  function by $I(.)$. We may regard the two  non-negative random variables $w_k$  and $v_k$ as two random weights  of the detection delay $(T-k)^+$ and  the event $I(T\geq k)$ such that the time of false alarm is greater than or equal to the change-point $k$, respectively. Here, $w_k, v_k \in \mathfrak{F}_{k-1}$ means that both weights $w_k$ and $v_k$  can be determined by the observation information before the time $k$ for $1\leq k\leq N$. Using the concept of the randomization probability of the change-point and the definition describing the average detection delay proposed by Moustakides (2008),  we can define a performance measure $\mathcal{J}_{M, N}(.)$ for every given weighted pair  $M=(W, V)$ to evaluate the detection performance of each sequential test $T\in \mathfrak{T}_N$ in the following
\begin{eqnarray}
\mathcal{J}_{M, N}(T)=\frac{\sum_{k=1}^{N+1}\textbf{E}_{k}(w_k(T-k)^+)}{\sum_{j=1}^{N+1}\textbf{E}_{0}(v_jI(T\geq j))}=\frac{\sum_{k=1}^{N}\textbf{E}_{k}(w_k(T-k)^+)}{\textbf{E}_{0}(\sum_{j=1}^{T}v_j)}.
\end{eqnarray}
Here, the second quality comes from $T\leq N+1$ and $\sum_{j=1}^{N+1}\textbf{E}_{0}(v_jI(T\geq j))= \textbf{E}_{0}(\sum_{j=1}^{T}v_j)$. As we only consider the detection delay after the change-point $\tau=k\geq 1$, the commonly-used detection delay $(T-k+1)^+$ is replaced by $(T-k)^+$ hereafter. Note that $W$ and $V$ may not be the randomization probability of the change-point.

 According to the definition of $ \mathcal{J}_{M, N}(T)$, the smaller $ \mathcal{J}_{M, N}(T)$, the better the detection performance of the test $T$ satisfying $\sum_{j=1}^{N+1}\textbf{E}_{0}(v_jI(T\geq j))\geq \gamma$ for some given positive constant $\gamma$.

\textbf{Remark 1.} The numerator  and  denominator of $\mathcal{J}_{M, N}(T)$ can be regarded as  a generalized out-of-control average run length ( ARL$_1$) and a generalized in-control ARL$_0$, respectively. Moreover, the measure $\mathcal{J}_{M, N}(.)$ can be considered as a generalization of the following Shiryaev's measure
\begin{eqnarray*}
\mathcal{J}_{S}(T)=\frac{\sum_{k=1}^{\infty}\rho_k\textbf{E}_{k}((T-k)^+)}{\sum_{j=1}^{\infty}\rho_j\textbf{E}_{0}(I(T\geq j))}=\textbf{E}(T-\tau|T\geq \tau).
\end{eqnarray*}
for $T\leq N+1$, where $\rho_k=\textbf{P}(\tau=k)$ for $k\geq 1$.

It is clear that taking various weighted pairs $M=(W,  V)$, we can get various measures $\mathcal{J}_{M, N}(.)$. Next we list six known measures  and two new measures in the following by taking the appropriate weighted pairs, $M_i=(W_i, V_i), 1\leq i\leq 8$.
\begin{eqnarray*}
\mathcal{J}_{M_1, N}(T)&=&\frac{\sum_{k=1}^{N+1}\rho_k\textbf{E}_{k}(T-k)^+}{\sum_{j=1}^{N+1}\rho_j\textbf{P}_{0}(T\geq j)},\,\,\,\,\,\,\,\, \mathcal{J}_{M_2, N}(T)=\frac{\textbf{E}_{1}(T-1)}{\textbf{P}_{0}(T\geq N+1)},\\
\mathcal{J}_{M_3, N}(T)&=&\frac{\textbf{E}_{1}(T-1)}{\sum_{j=1}^{N+1}\rho_j\textbf{P}_{0}(T\geq j)},\,\,\,\,\,\mathcal{J}_{M_4, N}(T)=\frac{\sum_{k=1}^{N}\textbf{E}_{k}((1-Y_{k-1})^+(T-k)^+)}{\textbf{E}_{0}(\sum_{j=1}^T(1-Y_{j-1})^+)},\\
\mathcal{J}_{M_5, N}(T)&=&\frac{r\textbf{E}_{1}(T-1)+\sum_{k=1}^{N}\textbf{E}_{k}((T-k)^+)}{r+\textbf{E}_{0}(T)},\\
\mathcal{J}_{M_6, N}(T)&=&\frac{\sum_{k=1}^{N}\textbf{E}_{k}((1-Y_{k-1})^+(T-k)^+)}{\textbf{E}_{0}(T)},\\
\mathcal{J}_{M_7, N}(T)&=&\frac{\textbf{E}_{1}((T-1))+\sum_{k=2}^{N+1}\textbf{E}_{k}(e^{X_{k-1}}(1+e^{X_{k-1}})^{-1}(T-k)^+)}{1+\textbf{E}_{0}(\sum_{k=2}^Te^{X_{k-1}}(1+e^{X_{k-1}})^{-1})},\\
\mathcal{J}_{M_8, N}(T)&=&\frac{\textbf{E}_{1}((T-1))+\sum_{k=2}^{N}\textbf{E}_{k}(\frac{1}{k-1}\sum_{j=1}^{k-1}e^{X_{j}}(T-k)^+)}{\textbf{E}_{0}(T)},
\end{eqnarray*}
where $Y_k=\max\{1, Y_{k-1}\}\Lambda_k$  for $1\leq k\leq N$, are the statistics of the CUSUM test with $Y_0=0$ (see  Moustakides 1986). For example, taking $W_1=V_1=\{\rho_k, 1\leq k\leq N+1\}$, $W_2=\{w_1=1, w_k=0, \, 2\leq k\leq N+1\}$, $V_2=\{v_j=0, \, 1\leq j\leq N, \, v_{N+1}=1\}$ and $W_4=V_4=\{w_j=v_j=(1-Y_{j-1})^+,\, 1\leq j\leq N+1\}$, we can get $\mathcal{J}_{M_i, N}(T)$ for $i=1, 2, 4.$  Since the in-control  ARL$_0$, $\textbf{E}_{0}(T)$, in $\mathcal{J}_{M_6, N}(T)$ is easier to be calculated than the generalized in-control ARL$_0$, $\textbf{E}_{0}(\sum_{j=1}^T(1-Y_{j-1})^+),$ in $\mathcal{J}_{M_4, N}(T)$, we often use the measure  $\mathcal{J}_{M_6, N}(T)$ to replace the measure $\mathcal{J}_{M_4, N}(T)$. Here,  both $e^{X_{k-1}}/(1+e^{X_{k-1}})$ and $\frac{1}{k-1}\sum_{j=1}^{k-1}e^{X_{j}}$ in the two new measures $\mathcal{J}_{M_7, N}(T)$ and $\mathcal{J}_{M_8, N}(T)$, can describe the changes of the observation values at change-point $k-1$ and the average of the changes of the observation values before the change-point $k\geq 2$, respectively.

Note that when we have an infinite independent observation sequences, the five measures above $\mathcal{J}_{M_i, \infty}$ for $ 1\leq i\leq 5$ and $N=\infty$, have been used by Shiryaev (1978, P. 193-200), Chow, Robbins and Siegmund (1971, P.108), Fris\'{e}n (2003),  Moustakides (1986) and Polunchenko and Tartakovsky (2010) to prove the optimality of the sequential tests, $T_S(c_1)$, $T_{SLR_1}(c_2)$, $T_{SLR_2}(c_3)$, $T_C(c_4)$  and $T^r_{SR}(c_5)$, respectively.

\subsection{Optimal sequential tests}

For a given weighted pair $M=(W, V)$, we first provide a definition of the optimization criterion of the sequential tests for  $N$ observations.

\textbf{Definition 1.}  \textit{ A sequential test $T^* \in \mathfrak{T}_N$ with $\textbf{E}_{0}(\sum_{k=1}^{T^*}v_k)\geq \gamma$ is  optimal under the measure $\mathcal{J}_{M, N}(T)$ if
\begin{eqnarray}
\inf_{T\in \mathfrak{T}_N,\,\,\,\textbf{E}_{0}(\sum_{j=1}^{T}v_j)\geq \gamma }\mathcal{J}_{M, N}(T)=\mathcal{J}_{M, N}(T^*)
\end{eqnarray}
where $\gamma$ satisfies $\textbf{E}_{0}(v_1)<\gamma <\textbf{E}_{0}(\sum_{j=1}^{N+1}v_j)$.}

To construct the optimal sequential test under the measure $\mathcal{J}_{M, N}(T)$ in (3) with a given weighted pair $M=(W, V)$, we need to present a series of nonnegative test statistics $\{Y_n, 0\leq n\leq N+1\}$,  as follows
\begin{eqnarray}
Y_n=\sum_{k=1}^nw_k\prod_{j=k}^n\Lambda^{(k)}_j
\end{eqnarray}
for $0\leq n\leq N+1$, where $Y_0=0$,  $Y_{N+1}=Y_N$, $W=\{w_k, 1\leq k\leq N+1\}$ and  $\Lambda^{(k)}_j$ satisfying (1).  It can be seen that the statistics $Y_n, 1\leq n\leq N,$ depend not only on  the likelihood ratio $\{\Lambda^{(k)}_j\}$ but also on the weight  of the detection delay $\{w_k\}$. Especially, if $\Lambda^{(k)}_j=\Lambda_j$ for $1\leq k\leq i\leq N$, that is, the post-change conditional densities (probabilities)  of the observation sequences do not depend on the change-point, then
\begin{eqnarray}
Y_n=\sum_{k=1}^nw_{k}\prod_{j=k}^n\Lambda_j=(Y_{n-1}+w_{n})\Lambda_n
\end{eqnarray}
for $1\leq n\leq N$.

\textbf{Remark 2.}  Even if (5) holds,  the test statistic sequence $\{Y_n, 0\leq n\leq N\}$ is  not necessarily a  Markov chain. For example,  let both the pre-change observation sequence $X_1, ...,X_{k-1}$  and the post-change observation sequence $X_k,... X_N$, be i.i.d., therefore, (5) holds, it is clear that the statistic $\{Y_n, 0\leq n\leq N\}$  is  not a Markov chain when we take  $w_1=1, w_n=\frac{1}{n-1}\sum_{j=1}^{n-1}e^{X_{j}}$ for $2\leq n\leq N$  in (5).

Motivated by Chow-Robbins-Siegmund's method of backward induction (1971, P.49), we present a  nonnegative random  dynamic  control limit $\{l_n(c), 0\leq n\leq N+1\}$ that is defined by the following recursive equations
\begin{eqnarray}
l_{N+1}(c)&=&0,\,\,\, l_N(c)=cv_{N+1} \nonumber \\
l_{n}(c)&=&cv_{n+1}+\textbf{E}_{0}\Big([l_{n+1}(c)-Y_{n+1}]^+|\mathfrak{F}_{n}\Big)
\end{eqnarray}
for $0\leq n\leq N-1$, where $c>0$ is a constant and $V=\{v_j, 1\leq j\leq N+1\}$. It is clear that $l_n(c)\geq cv_{n+1}$ and $ l_n(c) \in \mathfrak{F}_{n}$ for $0\leq n\leq N$.  The positive number $c$ can be regarded as an adjustment coefficient for the random dynamic control limit, as $l_n(c)$ is increasing on $c\geq 0$ with $l_n(0)=0$ and $\lim_{c\to \infty}l_n(c)=\infty$ for $v_{n+1}>0$.

Now, for a given weighted pair $M=(W, V)$,  we define a sequential test $T^*_{M}(c, N)$  by using the test statistics, $Y_n, 1\leq n\leq N+1,$ and the  control limits, $l_n(c), 1\leq n\leq N+1,$ as follows
\begin{eqnarray}
T^*_{M}(c, N)=\min\{ 1\leq n\leq N+1: \, Y_n\geq l_n(c)\}.
\end{eqnarray}
It is easy to check that $T^*_{M}(c, N)\in \mathfrak{T}_N$.

The following theorem shows that for any given performance measure $\mathcal{J}_{M, N}$ in (2), the sequential test $T^*_{M}(c, N)$ constructed above is optimal.

\textbf{Theorem 1.}  \textit{ Assume  that the ratio $\Lambda^{(k)}_j$ satisfies (1) for $1\leq k\leq N$ and $k\leq j\leq N$.  Let $\gamma$ be a positive number satisfying $\textbf{E}_0(v_1)<\gamma <\sum_{j=1}^{N+1}\textbf{E}_0(v_j)$. Then\\
(i) There exists a positive number $c_{\gamma}$
such that $T^*_{M}(c_{\gamma}, N)$ is optimal in the sense of (2) with $\textbf{E}_0(\sum_{j=1}^{T^*_{M}(c_{\gamma}, N)}v_j)=\gamma$; that is,
\begin{eqnarray}
\inf_{T\in \mathfrak{T}_N, \,\, \textbf{E}_{0}(\sum_{j=1}^{T}v_j)\geq \gamma }\mathcal{J}_{M, N}(T)=\mathcal{J}_{M, N}(T^*_{M}(c_{\gamma}, N)).
\end{eqnarray}
(ii) If $T\in \mathfrak{T}_N$ satisfies  $T\neq  T^*_{M}(c_{\gamma}, N)$, that is, $\textbf{P}_{0}( T\neq  T^*_{M}(c_{\gamma}, N))>0$  and  $\textbf{E}_{0}(\sum_{j=1}^{T}v_j)=\gamma$, then
\begin{eqnarray}
\mathcal{J}_{M, N}(T)>\mathcal{J}_{M, N}(T^*_{M}(c_{\gamma}, N)).
\end{eqnarray}
(iii) Moreover
\begin{eqnarray}
\mathcal{J}_{M, N}(T^*_{M}(c_{\gamma}, N))=c_{\gamma}\Big(1-\frac{\textbf{E}_{0}(v_1)}{\gamma}\Big)-\frac{\textbf{E}_{0}[l_1(c_{\gamma})-Y_1]^+}{\gamma}.
\end{eqnarray}}

Here, the random dynamic control limit $\{l_n(c), 0\leq n\leq N+1\}$ of the optimal test $T^*_{M}(c, N)$ can be called an optimal dynamic control limit.

It follows from (8) and (10) that  the  minimum value of the generalized out-of-control ARL$_1$  ( the numerator of the measure $\mathcal{J}_{M, N}(T)$ ) for all $T\in \mathfrak{T}_N$ can be calculated using the following formula
\begin{eqnarray}
&&\inf_{T\in \mathfrak{T}_N, \,\, \textbf{E}_{0}(\sum_{j=1}^{T}v_j)\geq \gamma }\sum_{k=1}^{N}\textbf{E}_{k}(w_k[T-k]^+)\nonumber\\
&=&\sum_{k=1}^{N}\textbf{E}_{k}(w_k[T^*_{M}(c_{\gamma}, N)-k]^+)\\
&=&c_{\gamma}(\gamma  -\textbf{E}_{0}(v_{1}))-\textbf{E}_{0}([l_{1}(c_{\gamma})-Y_{1}]^+).\nonumber
\end{eqnarray}

\textbf{Remark 3.}   Unless the test statistic $ \{Y_k, 0\leq n\leq N\}$ is a Markov chain, it is hard to prove the optimality of $T^*_{M}(c, N)$ under the measure $\mathcal{J}_{M, N}(T)$ by the optimal stopping method of Markov sequences proposed by Shiryaev (1978).

As an application of Theorem 1, we have the following corollary.

\textbf{Corollary 1.}  \textit{ The eight sequential tests $T^*_{M_i}(c, N), 1\leq i \leq 8$, defined in (7), which correspond to the eight weighted pairs $M_i, 1\leq i\leq 8$, are optimal under the measures $\mathcal{J}_{M_i, N}$ listed in Section 2.1 for $1\leq i\leq 8$, respectively. }

Note that the optimality of the two tests $T^*_{M_4}(c, N)$ and $T^*_{M_6}(c, N)$ with the optimal dynamic control limits $\{l^{(4)}_n(c), 0\leq n\leq N+1\}$ and $\{l^{(6)}_n(c), 0\leq n\leq N+1\}$, respectively, is not under Lorden's measure (see Lorden 1971,  Moustakides 1986) but under the corresponding measures $\mathcal{J}_{M_4, N}$ and $\mathcal{J}_{M_6, N}$, respectively.

\section{Optimal control limits}

It is clear that the optimal control limit $\{l_n(c), 0\leq n\leq N+1\}$ of the optimal sequential test $T^*_{M}(c, N)$ plays a key role in detecting changes in distribution. Since $\textbf{E}_{0}([l_{n+1}(c)-Y_{n+1}]^+|\mathfrak{F}_{n})$ and $v_{n+1}$ are measurable with respect to $\mathfrak{F}_{n}$, it follows that there are $2N+1$ non-negative functions $h_n=h_n(c, x_0, x_1, ..., x_n)$, $0\leq n\leq N-1,$  and $v_{n}=v_{n}(x_0, x_1, ...,x_{n-1}),$ $0\leq n\leq N-1,$ such that
\begin{eqnarray*}
h_n&=&h_n(c, x_0, x_1, ..., x_n)\\
&=&\textbf{E}_{0}([l_{n+1}(c)-Y_{n+1}]^+|X_n=x_n, X_{n-1}=x_{n-1}, ..., X_0=x_0)
\end{eqnarray*}
for $0\leq n\leq N-1$. Therefore, the optimal control limit $l_n(c)$ in (6) can be written as
\begin{eqnarray*}
l_n(c)=cv_{n+1}(x_0, X_1, ...,X_n) +h_n(c, x_0, X_1, ..., X_n)
\end{eqnarray*}
for $0\leq n\leq N$, where $X_0=x_0$ is a constant. It can be seen that the optimal control limit $\{l_n(c), 0\leq n\leq N+1\}$ of the optimal sequential test $T^*_{M}(c, N)$ is not easy to calculate for a general dependence observation sequence $\{X_n, \, 0\leq n\leq N\}$.

To reduce the number of observation variables on which the control limit $\{l_n(c), \, 0\leq n\leq N\}$ depends, we let the observation sequence   $\{X_n, \, 0\leq n\leq N\}$  be at most a $\textbf{p}$-order Markov chain, where $p=\max\{i,\,j\}$, $0\leq p\leq N$, that is, both the pre-change observations $X_1, ..., X_{k-1}$ and  the post-change observations $X_{k}, ..., X_{N}$ are $i$-order and $j$-order Markov chains with transition probability density functions $p_{\textbf{0}n}(x_{n}|x_{n-1},..., x_{n-i})$ and  $p^{(k)}_{\textbf{1}m}(x_{m}|x_{m-1},..., x_{m-j})$, respectively, which satisfy the following Markov  property
\begin{eqnarray*}
p_{\textbf{0}n}(x_{n}|x_{n-1},..., x_{n-i})&=&p_{\textbf{0}n}(x_{n}|x_{n-1},..., x_{n-i},...,x_0)\\
p_{\textbf{1}m}(x_{m}|x_{m-1},..., x_{m-j})&=&p_{\textbf{1}m}(x_{m}|x_{m-1},..., x_{m-j},...,x_0)\\
&=& p^{(k)}_{\textbf{1}m}(x_{m}|x_{m-1},..., x_{m-j},...,x_0)
\end{eqnarray*}
for $n\geq i$, $m\geq j$ and $1\leq k\leq m\leq N$. The last equation above means that the post-change conditional densities  of the observation sequences do not depend on the change-point. Here, a $\textbf{0}$-order Markov chain means that both the pre-change observations $X_1, ..., X_{k-1}$ and  the post-change observations $X_{k}, ..., X_{N}$ are mutually independent. When $p=N$, we consider that  at least one of the pre-change observations $X_1, ..., X_{k-1}$ and the post-change observations $X_{k}, ..., X_{N}$ is not a Markov chain of any order since we have only $N$ observations. In this case, the test statistics, $Y_0, Y_1,..., Y_N$, can be considered not to be a Markov chain of any order even if the post-change conditional densities  of the observation sequences  do not depend on the change-point.

The following theorem 2 shows that the optimal control limit $l_n(c) ( 0\leq n\leq N)$ depends on $Y_n$ and $p$ observation variables, if the observation sequence $\{X_n, \, 0\leq n\leq N\}$  is at most a $\textbf{p}$-order Markov chain.

\textbf{Theorem 2.}  \textit{ Let the observation sequences be at most a $\textbf{p}$-order Markov chain for $0\leq p\leq N$. Let $A_{n, p}:=\{X_{n},..., X_{n-p+1}\}$ and $A_{n, 0}:=\{X_{n},..., X_{0}\}$. Assume that the post-change conditional densities  of the observation sequences  do not depend on the change-point and the weighted pair $M=(W, V)$ satisfy $w_{n+1}=w_{n+1}(Y_{n}, A_{n, p_1})$ and $v_{n+1}=v_{n+1}(Y_{n},  A_{n, p_2})$ for  $0\leq n\leq N$, where $0\leq p_1, p_2\leq p$, $w_{n+1}=w_{n+1}(Y_{n})$ for $p_1=0$ and $v_{n+1}=v_{n+1}(Y_{n})$ for $p_2=0$. Then \\
(i) For $1\leq p\leq N$, the optimal control limit $\{l_n(c), \, 0\leq n\leq N\}$ can be written as
\begin{eqnarray*}
l_n(c)&=&cv_{n+1}(Y_{n}, A_{n, p_2})\\
&+&\textbf{E}_{0}\Big([l_{n+1}(c)-(Y_n+w_{n+1}(Y_{n}, A_{n, p_1}))\Lambda_{n+1}]^+|Y_n, A_{n, 0}\Big)
\end{eqnarray*}
for $0\leq n\leq p-1$ and
\begin{eqnarray*}
l_n(c)&=&cv_{n+1}(Y_{n}, A_{n, p_2})\\
&+&\textbf{E}_{0}\Big([l_{n+1}(c)-(Y_n+w_{n+1}(Y_{n}, A_{n, p_1}))\Lambda_{n+1}]^+|Y_n, A_{n, p}\Big)
\end{eqnarray*}
for $p\leq n\leq N$, where we will replace $X_{n-p_1+1}$ or $X_{n-p_2+1}$ with $X_0$ as long as $n-p_1+1<0$ or  $n-p_2+1<0$ respectively. \\
(ii) For $p=0$, we have
\begin{eqnarray*}
l_n(c)=cv_{n+1}(Y_n)+\textbf{E}_{0}\Big([l_{n+1}(c, Y_{n+1})-(Y_n+w_{n+1}(Y_n))\Lambda_{n+1}]^+|Y_n\Big)
\end{eqnarray*}
for $0\leq n\leq N$. }

Note that the optimal control limit $l_n(c)$  depends  not only  on $ A_{n, p}$ but also  on the test statistic  $Y_n$ for $1\leq n\leq N$. Can we find a control limit $\widetilde{l_n}(c)$ that has the same property as $l_n(c)$ but does not directly depend on the test statistic $Y_n$  for $1\leq n\leq N$ ?  To answer this question, we first give a definition of  an equivalent control limit.

\textbf{Definition 2.}  \textit{ Let the observation sequence $\{ \widetilde{l_n}(c),\, 1\leq n\leq N\}$ be a control limit of  a sequential test  $\widetilde{T} \in \mathfrak{T}_N$, where $\widetilde{T} =\min\{ 1\leq n\leq N+1: \, Y_n\geq \widetilde{l_n}(c)\}$. If $\widetilde{T}$ is equal to the optimal sequential test $T^*_{M}(c, N)$ ( a.s. $\textbf{P}_{0}$ ), then we call the control limit $\{ \widetilde{l_n}(c) \}$ an equivalent control limit of the optimal sequential test $T^*_{M}(c, N)$.}

The following theorem answers the above question.

\textbf{Theorem 3.}  \textit{ Let  the observation sequences  and the weighted pair $M=(W, V)$ satisfy the conditions of Theorem 2.  Let $a_{n, p}:=\{x_{n},..., x_{n-p+1}\}$ and $a_{n, 0}:=\{x_{n},..., x_{0}\}$. Assuming that $p_1=p_2=p$, $ y+w_{n+1}(y, a_{n, p})$ and $v_{n+1}(y, a_{n, p})$  are  continuous nondecreasing and non-increasing  on $y\geq 0$ respectively for given $a_{n, p}$, $0\leq n\leq N$. Then \\
(i) For $1\leq p\leq N$, there is an equivalent control limit $\widetilde{l_n}(c)$ of the optimal sequential test $T^*_{M}(c, N)$ which  does not depend directly on the statistic $Y_n$ for $1\leq n\leq N$ such that $\widetilde{l_n}(c)=y_n(c, A_{n, 0})$ for $0\leq n\leq p-1$ and $\widetilde{l_n}(c)=y_n(c, A_{n, p})$ for $p\leq n\leq N$, where the nonnegative functions $y_n= y_n(c, a_{n, 0})$ for $0\leq n\leq p-1$ and $y_n=y_n(c, a_{n, p})$ for $p\leq n\leq N$ satisfy the following equations
\begin{eqnarray*}
y_n&=&cv_{n+1}(y_n, a_{n, 0})\\
&+&\textbf{E}_{0}\Big([l_{n+1}(c)-(y_n+w_{n+1}(y_n, a_{n, 0}))\Lambda_{n+1}]^+|Y_n=y_n, A_{n, 0}=a_{n, 0}\Big)
\end{eqnarray*}
for $0\leq n\leq  p-1$ and
\begin{eqnarray*}
y_n&=&cv_{n+1}(y_n, a_{n, p})\\
&+&\textbf{E}_{0}\Big([l_{n+1}(c)-(y_n+w_{n+1}(y_n, a_{n, p}))\Lambda_{n+1}]^+|Y_n=y_n, A_{n, p}=a_{n, p}\Big)
\end{eqnarray*}
for $ p \leq n\leq N$. \\
(ii) Let $p=0$. There is a series of nonnegative non-random numbers, $y_n, 1\leq n\leq N,$ such that the equivalent control limit $\widetilde{l_n}(c)=y_n$ and $y_n$ satisfies
\begin{eqnarray}
     \,\,\,\,\,\,\,\,\,\,\,\,\,\,\,\, y_n=cv_{n+1}(y_n)+\textbf{E}_{0}\Big([l_{n+1}(c)-(y_n+w_{n+1}(y_n))\Lambda_{n+1}]^+|Y_n=y_n\Big)
\end{eqnarray}
for $1\leq n\leq N$.}

\textbf{Remark 4.}  By the similar method of proving Theorem 3, we can prove that the results of Theorem 3 are still true for $0\leq p_1, p_2\leq p\leq N$.

It is clear that the weighted pairs $M_i$ for $1\leq i\leq 6$, satisfy the conditions of Theorem 3. As an application of Theorem 3, we have the following corollary.

\textbf{Corollary 2.}  \textit{  Let the observation sequences be at most a $\textbf{p}$-order Markov chain for $0\leq p\leq N$ and the post-change conditional densities  of the observation sequences  do not depend on the change-point. Then,  the six optimal sequential tests $T^*_{M_i}(c, N)$ for $1\leq i\leq 6$, have equivalent control limits. Especially, when $p=0$, the equivalent control limits  consist of a series of dynamic non-random numbers. }

Since  none of the equivalent control limits of optimal sequential tests $T^*_{M_i}(c, N)$ for $1\leq i\leq 6$ are constants when $p=0$. This means that $T^*_{M_1}(c, N) \neq T_{S, N}(c_1)$, $T^*_{M_2}(c, N) \neq T_{SLR_1, N}(c_2)$, $T^*_{M_3}(c, N) \neq T_{SLR_2, N}(c_3)$,  $T^*_{M_4}(c, N) \neq T_{C, N}(c_4)$ and $T^*_{M_5}(c, N) \neq T^r_{SR, N}(c_5)$, since the control limits, $c_i, 1\leq i\leq 5,$ are constants.

Thus, from (ii) of Theorem 1, we can get the following corollary

\textbf{Corollary 3.}  \textit{  The optimal sequential tests $T^*_{M_i}(c, N)$ for $1\leq i \leq 5,$ are strictly superior to the tests $ T_{S, N}(c_1)$, $T_{SLR_1, N}(c_2)$, $T_{SLR_2, N}(c_3)$, $T_{C, N}(c_4)$ and  $T^r_{SR, N}(c_5)$ under the measures $\mathcal{J}_{M_i, N}$ for  $1\leq i\leq 5$, respectively,  when they all have the same  (generalized) in-control ARL$_0$.}

\textbf{Remark 5.}  Sections 4.1 and 4.2 illustrate that the CUSUM and Shiryaev-Roberts tests with appropriate dynamic control limits can be superior to the  CUSUM test  $T_{C, N}$ under Lorden's measure and the Shiryaev-Roberts test $T^r_{SR, N}$ under Pollak's measure (see Pollak 1985) respectively for finite independent observations. Thus, the reason why the optimal sequential tests mentioned in the Introduction, $T_S(c_1)$, $T_{SLR_1}(c_2)$, $T_{SLR_2}(c_3)$, $T_C(c_4)$ and $T^r_{SR}(c_5)$  for a sequence of infinite independent observations are no longer optimal for finite independent observation sequences, is that all of their control limits, $c_k, 1\leq k\leq 5$, are constants.

Next, we illustrate how to find an equivalent control limit by analyzing the optimal control limit of the optimal sequential test $T^*_{M_2}(c, N)$. Let the observation sequence   $\{X_k, \, 0\leq k\leq N\}$  be independent, that is, $p=0$. Take $W_2=\{w_1=1, \, w_k=0, \, 2\leq k \leq N+1\}$ and  $V_2=\{v_k=0, \, 1\leq k\leq N, \, v_{N+1}=1\}$, we know that $\{Y_n=\prod_{j=1}^n\Lambda_{j}, \, 1\leq n\leq N\}$ is a Markov chain and  $\Lambda_{n+1}=p_{\textbf{1}(n+1)}(X_{n+1})/p_{\textbf{0}(n+1)}(X_{n+1})$ and $Y_n$ are mutually independent for $1\leq n\leq N-1$. It follows from (6) and (ii) of Theorem 2 that the optimal control limit of $T^*_{M_2}(c, N)$ can be written as
\begin{eqnarray*}
l_{N+1}(c)&=&0,\,\,\, l_N(c)=c \nonumber \\
l_{N-1}(c, y)&=&\textbf{E}_{0}\Big([c-y\Lambda_{N}]^+|Y_{N-1}=y\Big)=\textbf{E}_{0}\Big([c-y\Lambda_{N}]^+\Big)\\
l_{n}(c, y)&=&\textbf{E}_{0}\Big([l_{n+1}(c, Y_{n+1})-Y_{n+1}]^+|Y_{n}=y\Big)\\
&=&\textbf{E}_{0}\Big([l_{n+1}(c, y\Lambda_{n+1})-y\Lambda_{n+1}]^+\Big)
\end{eqnarray*}
for $1\leq n\leq N-2.$ It is clear that the function $l_{N-1}(c, y)$ is strictly monotonically decreasing on $y\geq 0$. Hence, $l_{n}(c, y)$ is also strictly monotonically decreasing on $y\geq 0$ for  $1\leq n\leq N-2.$ This means that for each $k$ ($1\leq n\leq N-1$), there is a unique positive number $y_n$ such that $y_n=l_{n}(c, y_n)$ for $c>0$. Thus, $Y_n\geq y_n$ if and only if $Y_n\geq l_{n}(c, Y_n)$ for $1\leq n\leq N-1$. In other words, the equivalent control limits $\{\widetilde{l_n}(c),\, 1\leq n\leq N\}$ of the optimal sequential test $T^*_{M_2}(c, N)$  are a series of positive numbers $\{y_n, \,1\leq n\leq N\}$, where $y_N=c>0$ and $y_n$ satisfies $y_n=l_{n}(c, y_n)$ for $1\leq n\leq N-1$.

\section{ Comparison and analysis of simulation results}

Consider an observation sequence with $N=60$. Let the change time $\tau$ be unknown. By comparing the simulation results respectively in Sections 4.1 and 4.2, we illustrate that the CUSUM test $T_C$  and  the Shiryaev-Roberts test  $T^r_{SR}$  with a specially designed deterministic initial point $r$ for an exponential model,  are no longer optimal under Lorden's and Pollak's measures for $60$ finite independent observations, respectively. The detection performance (the generalized  out-of-control ARL$_1$)  of six sequential tests, $T^*_{M_5}(c, 60)$, $T^*_{M_6}(c, 60)$, $T_C$, $T_{E}$,  $T^{-1/60}_C$  and $T^{1/60}_C$ for the independent or dependent observation sequence,  are compared  in Sections 4.3 and 4.4 respectively,  where   $T_E$ denotes the EWMA ( the exponentially weighted moving average ) test introduced by Roberts (1959), which, like the CUSUM test $T_C$, is very popular in statistical process control (see  Han and Tsung, 2004; Saleh, \emph{et al}. 2015; Hosseini and Noorossana, 2018). Both $T^{-1/60}_C$  and $T^{1/60}_C$ are defined by replacing  the constant control limit of the CUSUM test $T_C$ with two straight lines, $c^-_k=c(1-k/60)$ and $c^+_k=c(1+k/60)$ for $1\leq k\leq 60$, respectively. All the numerical simulation results in this section were obtained using $\textbf{10}^5$  repetitions.

\subsection{ Comparison of simulation values of $\mathcal{J}_L(\min\{T, N+1\})$}

Let $\{X_k, 1\leq k\leq 60\}$ be an i.i.d observation sequence  with a pre-change normal  distribution of $N(0, 1)$ and a post-change normal distribution of $N(0.2, 1)$. That is, the likelihood ratio $\Lambda_k$ of the pre-change and post-change probability densities $p_0(x)$ and $p_1(x)$ can be written as $\Lambda_k=e^{0.2(X_k-0.1)}$ for $1\leq k\leq 60$. We will compare the performance of the two CUSUM tests $T_C(c, 60)$ and  $T_{DC}$ in detecting the mean shift from $\mu_0=0$ to $\mu_1=0.2$  under Lorden's measure $J_L(\min\{T, N+1\})$ with ARL$_0$=40, where $T_C(c_4, 60)=\min\{T_C(c_4), 61\}$ and
\begin{eqnarray*}
T_{DC}=\min\{ 1\leq k\leq N+1: \, Y_k\geq l_k\},
\end{eqnarray*}
with the following dynamic control limits
\begin{equation*}
   l_k=\begin{cases}
   2.53 &\mbox{if $1 \leq k \leq 40$ }\\
   2.53+0.506*(k-40) &\mbox{if $40<k \leq 60,$}
   \end{cases}
\end{equation*}
and $l_{61}=0$,  where  $Y_{61}=Y_{60}$, $Y_k, \, 0\leq k \leq 60,$ are the CUSUM test statistics, that is, $Y_0=0$ and $Y_k=\max\{1, Y_{k-1}\}\Lambda_k$ for $1\leq k\leq 60$. It can be calculated that $\textbf{E}_{0}(T_{DC})=40.02$.

Taking the constant control limit $c_4=2.6601$, we have $\textbf{E}_{0}(T_C(c_4, 60))=40.01.$ Note that
\begin{eqnarray*}
essup\{\textbf{E}_{k}((T_C(c_4, 60)-k)^+|\mathfrak{F}_{k-1})\}=\textbf{E}_{k}((T_C(c_4, 60)-k)^+|Y_{k-1}\leq 1)
\end{eqnarray*}
for $1\leq k\leq 60$.  Both the simulation values of the detection delay $\textbf{E}_{k}((T_C(c_4, 60)-k)^+|Y_{k-1}\leq 1)$ and $\textbf{E}_{k}((T_{DC}-k)^+|Y_{k-1}\leq 1)$  are decreasing for $k=1, 2, ..., 60$, that is, both can arrive the  maximum values  at change-point $k=1$. Since  both $\textbf{E}_{1}(T_{DC}-1)=22.951$ and $\textbf{E}_{1}(T_C(c, 60)-1)=23.425$ are the maximum values, it follows that
\begin{eqnarray*}
&&\mathcal{J}_L(T_{DC})=\max_{1\leq k\leq 60}\{\textbf{E}_{k}((T_{DC}-k)^+|Y_{k-1}\leq 1)\}=\textbf{E}_{1}(T_{DC}-1)\\
&<& \mathcal{J}_L(T_C(c, 60))=\max_{1\leq k\leq 60}\{\textbf{E}_{k}((T_C(c, 60)-k)^+|Y_{k-1}\leq 1)\}= \textbf{E}_{1}(T_C(c, 60)-1).
\end{eqnarray*}
This means that the CUSUM chart $T_C$ is not optimal under Lorden's measure $\mathcal{J}_L(\min\{T, N+1\})$ restricted in  $60$ i.i.d. observation sequence.

\subsection{ Comparison of simulation values of $ \mathcal{J}_P(\min\{T, N+1\})$}

Let  $\{X_k, 1\leq k\leq 60\}$ be an i.i.d observation sequence  with a pre-change exponential density of $f_0(x)=e^{-x}I(x\geq 0)$ and a post-change exponential density of $f_1(x)=2e^{-2x}I(x\geq 0)$. The likelihood ratio is $\Lambda_k=2e^{-X_k}$ for $1\leq k\leq 60$. Polunchenko and Tartakovsky (2010) have proved that the control chart $T^r_{SR}(c)$ with a specially designed deterministic initial point $r$ for
an exponential model is optimal under Pollak's measure $\mathcal{J}_P(T)$ for $1<\gamma <2.2188$. Let $T^r_{SR}(c_5, 60)=\min\{T^r_{SR}(c_5), 61\}$. Taking $c_5=1.6645$ and $r=\sqrt{2.6645}-1$, we have ARL$_0=\textbf{E}_{0}(T^r_{SR}(c_5, 60))=2$.  It follows from $\mathcal{J}_P(\min\{T, N+1\})=\max_{1\leq k\leq 60}\{\textbf{E}_{k}(T-k)^+/\textbf{P}_{0}(T \geq k)\}$ that
\begin{eqnarray*}
\mathcal{J}_P(T^r_{SR}(c_5, 60))=\textbf{E}_{1}(T^r_{SR}(c_5, 60)-1)=1.3165
\end{eqnarray*}
However, if we define a sequential test as $ T^r_{SR}(\{l_k\}, 60)$ with dynamic control limit $l_k$
\begin{equation*}
   l_k=\begin{cases}
   1.238+0.1238k &\mbox{if $1 \leq k \leq 10$ }\\
   0 &\mbox{if $10<k \leq 60$},
   \end{cases}
\end{equation*}
we can obtain
\begin{eqnarray*}
\mathcal{J}_P(T^r_{SR}(\{l_k\}, 60))=\textbf{E}_{1}(T^r_{SR}(\{l_k\}, 60)-1)=1.2743
\end{eqnarray*}
with ARL$_0=\textbf{E}_{0}( T^r_{SR}(\{l_k\}, 60))=2.0012 $. Thus
\begin{eqnarray*}
\mathcal{J}_P(T^r_{SR}(\{l_k\}, 60)) < \mathcal{J}_P(T^r_{SR}(c_5, 60)).
\end{eqnarray*}
This means that the control chart $T^r_{SR}(c_5)$ is not optimal under Pollak's measure $\mathcal{J}_P(\min\{T, N+1\})$ restricted in  $60$ i.i.d. observations.

\subsection{ Comparison of the generalized  out-of-control ARL$_1$ for independent observations}

Let $\{X_k, 1\leq k\leq 60\}$ be an i.i.d. observation sequence  with a pre-change normal  distribution of $N(0, 1)$ and a post-change normal distribution of $N(1, 1)$.  The likelihood ratio is $\Lambda_k=e^{X_k-1/2}$ for $1\leq k\leq 60$. Let $T^*_5=T^*_{M_5}(c, 60)$, $T^*_6=T^*_{M_6}(c, 60)$ and let the smoothing parameter in the statistics of the EWMA test $T_E$ be 0.1. By Corollary 2, we know that the equivalent control limits of the optimal sequential tests $T^*_5$ and $T^*_6$ consist of a series of non-random positive numbers. Fig. 1 shows the constant control limit of $T_C$ (black dots) and the equivalent dynamic control limit of $T^*_6$ (white dots).

\begin{figure}[H]
  \centering
  \includegraphics[width=12cm]{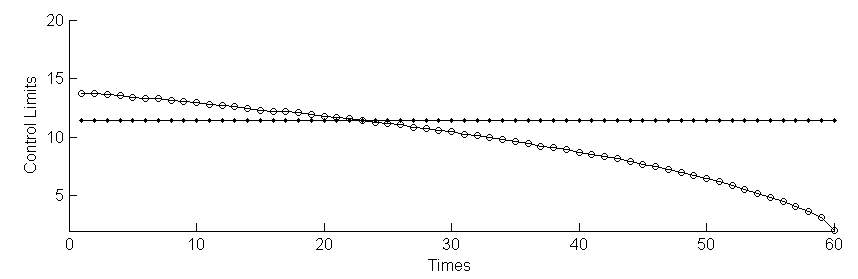}\\
  \caption{Control limits for $T_C$ and $T^*_6$ with ARL$_0\approx $ 40 }
  \end{figure}

We use two generalized out-of-control ARL$_1$s,  GARL$_5$ and GARL$_6$,  to evaluate the detection performance of the sequential tests, where
\begin{eqnarray*}
 GARL_5(T)&=&\textbf{E}_{0}(T)\mathcal{J}_{M_5, N}(T)=\sum_{k=1}^{N}\textbf{E}_{k}((T-k)^+)\\
 GARL_6(T)&=&\textbf{E}_{0}(T)\mathcal{J}_{M_6, N}(T)=\sum_{k=1}^{N}\textbf{E}_{k}((1-Y_{k-1})^+(T-k)^+)
\end{eqnarray*}
where $r=0$ in $\mathcal{J}_{M_5, N}(T)$. Obviously, for any two sequential tests $T', T\in \mathfrak{T}_N$ with $\textbf{E}_{0}(T')=\textbf{E}_{0}(T)$, we have  $GARL_j(T')\geq GARL_j(T)$ if and only if $\mathcal{J}_{M_j, N}(T')\geq \mathcal{J}_{M_j, N}(T)$ for $j=5, 6$.

The simulation results of GARL$_5$ and GARL$_6$ for the six tests $ T^*_5$, $T^*_7$, $T_C$, $T_{E}$, $T^{-1/60}_C$  and $T^{1/60}_C$ with the same ARL$_0\approx 20, 40, 50$, are listed in Table 1, where the values of ARL$_0$, the constant control limits of $T_C$ and $T_{E}$, and the adjustment coefficients of  $ T^*_5$, $T^*_6$, $T^{-1/60}_C$  and $T^{1/60}_C$  are listed in parentheses. Table 1 shows that both $T^*_5$ and $T^*_6$ have the best detection performance; that is, $T^*_5$  and $T^*_6$  have the smallest GARL$_5$ and GARL$_6$ (in bold) respectively  in the six tests with the same ARL$_0 \approx 20, 40, 50$.   This is consistent with the result of Corollary 3: tests  $T^*_5$  and $T^*_6$ are optimal  under measures $\mathcal{J}_{M_5, N}(T)$ and $\mathcal{J}_{M_6, N}(T)$ respectively.

\footnotesize
\setlength{\tabcolsep}{2pt}
\begin{tabular}{l|l|cccccc}
\multicolumn{8}{l}{\textbf{Table 1.}}Simulation values of GARL$_5$ and GARL$_6$ with the same ARL$_0$ for indepen-\\
\multicolumn{2}{l}{}dent observations\\
\hline\hline
\textbf{ARL$_0$}&  &  & &\textbf{Sequential Tests}  &  & & \\
& &\textbf{$T^*_5$}&\textbf{$T^*_6$}&\textbf{$T_C$}&\textbf{$T_{E}$}&\textbf{$T^{-1/60}_C$}&\textbf{$T^{1/60}_C$}\\
\hline\hline
& GARL$_5$&\textbf{42.10}&44.75&45.13&48.02&46.50&47.57 \\
& GARL$_6$&19.62&\textbf{17.59}&18.97&19.98 &19.28&19.34\\
\cline{2-8}
&c&(0.12216) &(1.3011)&(4.4823)&(1.2250)&(6.3900)&(3.629)\\
\raisebox{4.0ex}[0pt]{$20$}&ARL$_0$&(20.01)&(20.06)&(20.07)&(20.08)&(20.08)&(20.07)\\
\hline
&GARL$_5$&\textbf{139.18}&145.65&148.07&164.28&148.76&155.80 \\
&GARL$_6$& 55.17&\textbf{49.26}&54.44&59.97&54.96&55.99\\
\cline{2-8}
&c&(5.5996)&(2.0251)&(11.4423)&(1.4064)&(22.1500)&(8.7815)\\
\raisebox{4.0ex}[0pt]{$40$}&ARL$_0$&(40.02)&(40.06)&(40.06)&(40.04)&(40.01)&(40.02)\\
\hline
&GARL$_5$&\textbf{229.26}&232.52&240.52&273.29&238.82&248.57 \\
& GARL$_6$&84.27&\textbf{80.95}&83.45& 95.52 &83.85&85.63\\
\cline{2-8}
&c&(0.2656)&(2.9518)&(22.8821)&(1.5269)&(52.2500)&(17.2478)\\
\raisebox{4.0ex}[0pt]{$50$}&ARL$_0$&(50.02)&(50.05)&(50.04)&(50.08)&(50.00)&(50.05)\\
\hline\hline
\end{tabular}
\normalsize
\subsection{ Comparison of the generalized out-of-control ARL$_1$ for a Markov observation sequence}

Let $\{X_k, 1\leq k\leq 60\}$ be a dependent observation sequence satisfying
\begin{eqnarray*}
X_k= \displaystyle\left\{\begin{array}{clcc} \rho_0 X_{k-1} +\varepsilon_k
&\hspace{0.4cm}&{\hbox{if}}&
1\leq k\leq \tau, \hbox{} \\
\displaystyle  \rho_1 X_{k-1} +\varepsilon_k &\hspace{0.4cm}&{\hbox{if}}& k\geq \tau, \hbox{}
\end{array}\right.
\end{eqnarray*}
where $X_0=0$, $\{\varepsilon_k,\, 1\leq k\leq 60\}$ is i.i.d with a normal distribution, i.e.,   $\varepsilon_k \sim N(0,1)$ for $1\leq k\leq 60$, $\rho_0=0.5$ and $\rho_1=0.1$. That is, the correlation coefficient changes from $0.5$ to $0.1$. Obviously, $\{X_k, 1\leq k\leq 60\}$ is a Markov chain. The pre-change transition probability density $p_0(x, y)$, post-change transition probability density $p_1(x, y)$, and the likelihood ratio $\Lambda_k$ can be written respectively as
\begin{eqnarray*}
p_0(x, y)&=&\frac{1}{\sqrt{2\pi}}e^{-\frac{(y-\rho_0x)^2}{2}}, \,\,\,\,\,\,\,\,\,\, p_1(x, y)=\frac{1}{\sqrt{2\pi}}e^{-\frac{(y-\rho_1x)^2}{2}}\\ \Lambda_k&=&\frac{p_1(X_{k-1}, X_k)}{p_0(X_{k-1}, X_k)}=\exp\{[(\rho_1-\rho_0)X_{k-1}][X_k-(\rho_1+\rho_0)X_{k-1}/2]\}.
\end{eqnarray*}
It can be seen that the changes in the variance and covariance of $X_k$ and $X_{k-1}$ occur after the change-point $\tau=k$.  Here, the change-point is unknown.

As $\{X_k, 1\leq k\leq 60\}$ is a $\textbf{1}$-order  Markov chain, it follows from (i) of Theorem 3 that we need to calculate the equivalent control limits $\widetilde{l_k}=y_k(c, X_k)$ for $1\leq k\leq 59$ to get the corresponding optimal tests $ T^*_5$ and $T^*_6$ respectively.

We also use the two generalized out-of-control ARL$_1$s, GARL$_5$ and GARL$_6$, to evaluate the detection performance of the six sequential tests $ T^*_5$, $T^*_6$, $T_C$, $T_{E}$ with the smoothing parameter $0.1$, $T^{-1/60}_C$  and $T^{1/60}_C$. The simulation results of GARL$_5$ and GARL$_6$ for the six tests  with the same ARL$_0$=$ 20, 40$ and $50$, are listed in Table 2. The ARL$_0$ values, the constant control limits of $T_C$ and $T_{E}$, and the  adjustment coefficients of  $ T^*_5$, $T^*_6$, $T^{-1/60}_C$  and $T^{1/60}_C$  are listed in parentheses. Table 2 shows that tests $T^*_5$ and $T^*_6$ have the best detection performance; that is, $T^*_5$  and $T^*_6$  have the smallest GARL$_5$ and GARL$_6$ values (in bold) respectively of the six tests with the same ARL$_0 \approx 20, 40, 50$.  This is consistent with the result of Corollary 3: sequential tests $T^*_5$  and $T^*_6$ are optimal  under measures $\mathcal{J}_{M_5, N}(T)$ and $\mathcal{J}_{M_6, N}(T)$ respectively.

Note that though the monitoring performances of both $T^*_5$ and $T^*_6$ are better than all $T_C, T_E, T^{-1/60}_C$  and $T^{1/60}_C$ respectively under the measure $\mathcal{J}_{M_5, N}(T)$ and $\mathcal{J}_{M_6, N}(T)$, the constant control limits of $T_C, T_E, T^{-1/60}_C$ are easier to determine than that of $T^*_5$ and $T^*_6$.
\\
\\
\footnotesize
\setlength{\tabcolsep}{3pt}
\begin{tabular}{l|c|cccccc}
\multicolumn{8}{l}{\textbf{Table 2.}}Simulation values of GARL$_5$ and GARL$_6$ with the same ARL$_0$ for Markov\\
\multicolumn{3}{l}{}observation sequence\\
\hline\hline
& \textbf{}     &   & \textbf{Sequential Tests} &  &  & \\
\textbf{ARL$_0$} & \textbf{ARL$_1$} &\textbf{$T^*_5$}&\textbf{$T^*_6$}&\textbf{$T_C$}&\textbf{$T_{E}$}&\textbf{$T^{-1/60}_C$} &\textbf{$T^{1/60}_C$}\\
\hline\hline
&GARL$_5$&\textbf{115.43}&135.25&139.64&545.85&130.92&156.09 \\
&GARL$_6$&23.26&\textbf{21.55}&22.04&70.90 &22.72&23.09\\
\cline{2-8}
&c&(12.016)&(2.075)&(2.3482)&(0.7730)&(3.4500)&(1.8901)\\
\raisebox{4.0ex}[0pt]{$20$}&ARL$_0$&(20.05)&(20.14)&(19.97)&(20.06)&(20.01)&(20.09)\\
\hline
&GARL$_5$&\textbf{409.76}&467.17&474.64&665.57&450.68&1490.42 \\
&GARL$_6$&59.80&\textbf{57.86}&59.71& 151.81 &60.60&60.30\\
\cline{2-8}
&c&(22.8550)&(3.865)&(4.7828)&(0.933)&(10.3500)&(3.478)\\
\raisebox{4.0ex}[0pt]{$40$}&ARL$_0$&(40.72)&(40.84)&(40.76)&(40.03)&(40.02)&(40.03)\\
\hline
&GARL$_5$&\textbf{638.15}&688.52&705.62&1586.81&722.63&758.57 \\
&GARL$_6$&84.15&\textbf{80.42}&83.32& 181.89 &87.25&87.57\\
\cline{2-8}
&c&(32.89)&(5.575)&(7.528)&(1.0229)&(23.15)&(5.667)\\
\raisebox{4.0ex}[0pt]{$50$}&ARL$_0$&(49.77)&(49.26)&(49.28)&(49.99)&(49.94)&(50.04)\\
\hline\hline
\end{tabular}
\normalsize
\\
\section{Concluding remarks}

By presenting the generalized Shiryaev's measures of detection delay $\mathcal{J}_{M, N}(.)$, the statistic $Y_n, 0\leq n\leq N+1,$  the control limit $l_n(c), 0\leq n\leq N+1,$  and the sequential test  $T^*_{M}(c, N)$, we obtain the following main results. (i) For different measures $\mathcal{J}_{M, N}(.)$ of detection delay, we can construct different optimal sequential tests $T^*_{M}(c, N)$ under the corresponding measures for a general finite observation sequence. (ii) A formula is presented to calculate the value of the generalized out-of-control ARL$_1$ for every optimal test $T^*_{M}(c, N)$ which is the  minimum value of the generalized out-of-control ARL$_1$ of all test $T\in \mathfrak{T}_N$. (iii) When the post-change conditional densities (probabilities) of the observation sequences  do not depend on the change-point, there is an equivalent control limit that does not depend directly on the statistic of the optimal test $T^*_{M}(c, N)$ for $\textbf{p}$-order Markov chain. Specifically, the equivalent control limit can consist of a series of nonnegative non-random numbers when the observations are mutually independent.

In this paper, both the pre-change and post-change joint probability densities are assumed to be known. In fact, we usually do not know the post-change joint probability density before it is detected. But the potential
change domain (including the size and form of the
boundary) and its probability may be determined by
engineering knowledge and practical
experience. In other words, though the actual  post-change joint probability density  $p(\theta, k):=p_{\theta, k}(x_0, x_{1},...,
x_k,..., x_N)$ is unknown, that is, the parameter $\theta $ is unknown at the change time $k$,  we may assume that there is a known probability distribution  $Q_k(.)$ for the known parameter set $\Theta_{k}$ such that the  probability of the  post-change joint probability
density at change-point $k$ being $p_{\theta, k}$  is $dQ_k(\theta)$ for $1\leq k\leq N$, where  $p_{\theta, k}\neq p_{\theta', k}$ if and only if $\theta\neq \theta'$.  If we have no prior knowledge of the possible parameter $\theta$ (corresponding to a possible post-change probability density $p_{\theta, k}$) at the change-point $k$, it is natural to assume that  the probability distribution $Q_k$
may be an equal probability distribution or uniform
distribution on $\Theta_k$, that is, $Q_k(\theta=\theta_i)=1/m$ ($1\leq i\leq m<\infty$) for $\theta_i \in \Theta_k$ or $dQ(\theta)/d\theta=1/M(\Theta)$, where $dQ/d\theta$ denotes the probability density and $M(\Theta)$ is the measure (length, area, volume, etc.) of the bounded set $\Theta$. Note that the parameter $\theta$ may not be
the characteristic numbers (the mean, variance,
etc.) of the probability distribution.
Hence, we can define a new joint probability density
\begin{equation*}
p_{k}:=p_{k}(x_0, x_{1},...,
x_k,..., x_N)
\end{equation*}
 in the following
\begin{eqnarray*}
p_{k}(x_{0}, x_1, ..., x_N)=\int_{\Theta_{k}}p_{\theta, k}(x_{0}, x_1, ..., x_N)dQ_k(\theta)
\end{eqnarray*}
for $1\leq k\leq N$. The density function $p_{k}$ can be considered as a known post-change joint probability density at the change-point $k$, $1\leq k\leq N$.

\newpage
\textbf{APPENDIX : PROOFS OF THEOREMS}
\setcounter{equation}{0}
\renewcommand\theequation{A. \arabic{equation}}

{\bf Proof of Theorem 1. }  Since $(T-k)^+=\sum_{m=k+1}^{N+1}I(T\geq m)$ for $T\in \mathfrak{T}_N$, $I(T\geq m)\in \mathfrak{F}_{m-1}$ and the post-change joint probability density $p_{k}(x_0, x_1, ..., x_n)$ for the change-point $k$ ($1\leq k \leq N$) can be written as
\begin{eqnarray*}
p_{k}(x_0, x_1, ..., x_n)=q(x_0)\prod_{j=1}^{(k-1)\wedge n}p_{\textbf{0}j}(x_j|x_{j-1}, ..., x_0)\prod_{j=k}^np^{(k)}_{\textbf{1}j}(x_j|x_{j-1}, ..., x_0)
\end{eqnarray*}
for $1\leq n \leq N$, where $q(x_0)$ be the probability density (or probability) of $X_0$ at initial time $k=0$,  $(k-1) \wedge n$ denotes $\min\{k-1, n\}$, $\prod_{j=1}^{(k-1) \wedge n}=1$ for $k=1$ and $\prod_{j=k}^n=1$ for $n< k$,  it follows that
\begin{eqnarray}
&&\sum_{k=1}^{N}\textbf{E}_{k}(w_k(T-k)^+)\\
&=&\textbf{E}_{0}\Big(\sum_{k=1}^N\sum_{m=k+1}^{N+1}w_kI(T\geq m)\prod_{j=k}^{m-1}\Lambda^{(k)}_j\Big)\nonumber \\
&=&\textbf{E}_{0}\Big(\sum_{m=1}^{N}Y_{m}I(T\geq m+1)\Big)=\textbf{E}_{0}\Big(\sum_{m=1}^{T}Y_{m-1}\Big) \nonumber
\end{eqnarray}
for all $T\in \mathfrak{T}_N$. This equality means that the generalized out-of-control ARL$_1$ (the numerator of the measure $\mathcal{J}_{M, N}(T)$) is equal to the generalized in-control ARL$_0$, in which the weight $\{v_m\}$ is replaced by the statistic $\{Y_{m-1}\}$. Let $T^*=T^*(c)=T^*_M(c, N)$ and
\begin{eqnarray}
 \xi_n=\sum_{k=1}^n(Y_{k-1}-cv_k),
\end{eqnarray}
where $c>0$. We will divide three steps to complete the proof of Theorem 1.

\textbf{Step I.} Show that
\begin{eqnarray}
\textbf{E}_{0}(\xi_T) \geq \textbf{E}_{0}(\xi_{T^*})
\end{eqnarray}
for all $T\in \mathfrak{T}_N$ and the strict inequality of (A.3) holds for all $T\in \mathfrak{T}_N$ with $T\neq T^*$.

To prove (A.3), by Lemma 3.2 in Chow, Robbins and Siegmund (1971), we only need to prove the following two inequalities:
\begin{eqnarray}
\textbf{E}_{\infty}(\xi_{T^*}|\mathfrak{F}_{n}) \leq \xi_{n}\,\,\,\, \text{ on } \,\, \{ T^*> n\}
\end{eqnarray}
and
\begin{eqnarray}
\textbf{E}_{\infty}(\xi_{T}|\mathfrak{F}_{n}) \geq \xi_{n}\,\,\,\, \text{ on } \,\, \{T^*=n, \, T > n\}
\end{eqnarray}
for each $n\geq 1$.

Let $B_{m, n+1}(N)=\{Y_k < l_{k}(c), n+1\leq k\leq m\}$ for $n+1\leq m\leq N$. As in the proof of Theorem 1 in Han, Tsung and Xian (2017), we can  verify that
\begin{eqnarray}
l_{n}(c)=cv_{n+1}+\textbf{E}_{0}\Big(\sum_{m=n+1}^NB_{m, n+1}(N)[cv_{m+1}-Y_{m}]|\mathfrak{F}_{n}\Big)
\end{eqnarray}
and
\begin{eqnarray}
&&\textbf{E}_{0}\Big(\sum_{m=n+1}^NI(T>m)[cv_{m+1}-Y_m]|\mathfrak{F}_{n}\Big)\\
&\leq&(l_n(c)-cv_{n+1})I(T>n)\nonumber
\end{eqnarray}
for $0\leq n\leq N$ and $T\in \mathfrak{T}_N$. Note that here, $\{v_k, 1\leq k\leq N+1\}$ can be a series of non-negative random variables, but it is a series of positive numbers
$\{\rho_k, 1\leq k\leq N+1\}$ in the proof of Theorem 1 in Han, Tsung, and Xian (2017).

In fact, by the definition of $l_N(c)$ and $l_{N-1}(c)$  we know that (A.6) holds for $n=N-1$ and $n=N$, where $\sum_{m=n+1}^N=0$ for $n=N$.  Assume that (A.6) holds for $n\leq N-1$. Then, by the
definition of $l_{n-1}(c)$ and the assumption (A.6) for $l_{n}(c)$, we have
\begin{eqnarray*}
l_{n-1}(c)&=&cv_{n}+\textbf{E}_{0}\Big(B_{n, n}(N)[l_{n}(c)-Y_{n}]|\mathfrak{F}_{n-1}\Big)  \\
&=&cv_{n}+\textbf{E}_{0}\Big(B_{n, n}(N)[cv_{n+1}-Y_n]|\mathfrak{F}_{n-1}\Big)\\
  && +\textbf{E}_{0}\Big(\sum_{m=n+1}^NB_{n, n}(N)B_{m, n+1}(N)[cv_{m+1}-Y_{m}]|\mathfrak{F}_{n-1}\Big)\\
&=&  cv_{n}+\textbf{E}_{0}\Big(\sum_{m=n}^NB_{m, n}(N)[cv_{m+1}-Y_{m}]|\mathfrak{F}_{n-1}\Big).
\end{eqnarray*}
By mathematical induction, (A.6) holds  for $0\leq n\leq N$.  Furthermore, by (A.6) and the definition of $l_n(c)$ we have
\begin{eqnarray}
&&\textbf{E}_{0}\Big(\sum_{m=n+1}^NB_{m, n+1}(N)[cv_{m+1}-Y_{m}]|\mathfrak{F}_{n}\Big)\\
&&=l_{n}(c)-cv_{n+1}=\textbf{E}_{0}\Big([l_{n+1}(c)-Y_{n+1}]^+|\mathfrak{F}_{n}\Big)\nonumber
\end{eqnarray}
for $0\leq n\leq N$.

Obviously, (A.7) holds for $n=N$. For $n=N-1$, we have
\begin{eqnarray*}
&&\textbf{E}_{0}\Big(I(T>N)[cv_{N+1}-Y_N]|\mathfrak{F}_{N-1}\Big)\\
&\leq & \textbf{E}_{0}\Big(I(T>N)I(B_{N,N}(N))[cv_{N+1}-Y_N]|\mathfrak{F}_{N-1}\Big)\\
&\leq & I(T>N-1) \textbf{E}_{0}\Big(I(B_{N,N}(N))[l_N(c)-Y_N]|\mathfrak{F}_{N-1}\Big)\\
&=& (l_{N-1}(c)- cv_N)I(T>N-1),
\end{eqnarray*}
where the last equality comes from the definition of $l_{N-1}(c)$; that is, (A.7) holds for $n=N-1$. Assume that (A.7) holds for $n\leq N-1$. It follows that
\begin{eqnarray*}
&&\textbf{E}_{0}\Big(\sum_{m=n}^NI(T>m)[cv_{m+1}-Y_m]|\mathfrak{F}_{n-1}\Big)\\
&=&\textbf{E}_{0}(I(T>n)[cv_{n+1}-Y_n]|\mathfrak{F}_{n-1})\\
  &&+\textbf{E}_{0}\Big(\textbf{E}_{0}\Big(\sum_{m=n+1}^NI(T>m)[cv_{m+1}-Y_m]|\mathfrak{F}_{n}\Big)|\mathfrak{F}_{n-1}\Big)\\
&\leq &\textbf{E}_{0}(I(T>n)[cv_{n+1}-Y_n]|\mathfrak{F}_{n-1})+\textbf{E}_{0}\Big([ l_n(c)-cv_{n+1}]I(T>n)|\mathfrak{F}_{n-1}\Big)\\
&=&\textbf{E}_{0}\Big([l_n(c)-Y_n]I(T>n)|\mathfrak{F}_{n-1}\Big)\\
&\leq &\textbf{E}_{0}\Big(I(B_{n,n}(N))[l_n(c)-Y_n]I(T>n)|\mathfrak{F}_{n-1}\Big)\\
&\leq &I(T>n-1)\textbf{E}_{0}\Big(I(B_{n,n}(N))[l_n(c)-Y_n]|\mathfrak{F}_{n-1}\Big)\\
&=& I(T>n-1)(l_{n-1}(c)-cv_n),
\end{eqnarray*}
where the last equality comes from the definition of $l_{n-1}(c)$.  By mathematical induction, we know that (A.7) holds for all $0\leq n\leq N$.

From (A.6) and the definition of $\xi_{n}$ in (A.2), it follows that
\begin{eqnarray}
&&I(T^*> n)\textbf{E}_{0}(\xi_{T^*}-\xi_{n})|\mathfrak{F}_{n})\\
&=&I(T^*> n)\sum_{m=n}^{N}\textbf{E}_{0}(I(T^*> m)[\xi_{m+1}-\xi_{m}])|\mathfrak{F}_{n})\nonumber \\
&=&I(T^*> n)[Y_n-cv_{n+1}+\sum_{m=n+1}^{N}\textbf{E}_{0}(I(B_{m,n+1})(Y_m-cv_{m+1})]|\mathfrak{F}_{n}) \nonumber\\
&=&I(T^*> n)(Y_n-l_n(c))<0 \nonumber
\end{eqnarray}
for $1\leq n\leq N$. The last inequality comes from the definition of $T^*$. This means that (A.4) holds for $1\leq n\leq N$.

By (A.7), we have
\begin{eqnarray}
&&I(T^*= n)I(T>n)\textbf{E}_{0}(\xi_{T}-\xi_{n})|\mathfrak{F}_{n})\\
&=&I(T^*=n)\sum_{m=n}^{N}\textbf{E}_{0}(I(T>m)[\xi_{m+1}-\xi_{m}])|\mathfrak{F}_{n})\nonumber\\
&=&I(T^*=n)[I(T>n)(Y_n-cv_{n+1})\nonumber \\
&+&\sum_{m=n+1}^{N}\textbf{E}_{0}\Big(I(T>m)[Y_m-cv_{m+1}]|\mathfrak{F}_{n}\Big)]\nonumber\\
&\geq &I(T^*=n)[I(T>n)(Y_n-cv_{n+1})+( cv_{n+1}-l_n(c))I(T>n)]\nonumber\\
&= &I(T^*=n)I(T>n)[Y_n-l_n(c)]\geq 0.\nonumber
\end{eqnarray}
That is,  (A.5) holds for $1\leq n\leq N$. By (A.4) and (A.5), we know that  the inequality in (A.3) holds for all $T\in \mathfrak{T}_N$. Furthermore, from (A.9) and (A.10), it follows that the strict inequality in (A.3) holds for all $T\in \mathfrak{T}_N$ with $T\neq T^*$.

\textbf{Step II.} Show that there is positive number $c_{\gamma}$ such that
\begin{eqnarray}
\mathcal{J}_{M, N}(T^*(c_{\gamma}))=c_{\gamma}\Big(1-\frac{\textbf{E}_{0}(v_1)}{\gamma}\Big)-\frac{\textbf{E}_{0}[l_1(c_{\gamma})-Y_1]^+}{\gamma}. \nonumber
\end{eqnarray}
As  $\textbf{E}_0(v_1)<\gamma
<\sum_{k=1}^{N+1}\textbf{E}_0(v_k)$, it follows that there is at least a  $k\geq 2$ such that $\textbf{E}_0(v_k)>0$. Let $k^*=\max\{2\leq k\leq N+1: \,\,\textbf{E}_0(v_k)>0\}$, we have
\begin{eqnarray*}
\textbf{E}_{0}(\sum_{k=1}^{T^*}v_k)=\sum_{k=1}^{k^*}\textbf{E}_{0}(v_kI(T^*\geq k))=\textbf{E}_0(v_1)+\sum_{k=2}^{k^*}\textbf{E}_{0}(v_kI(T^*\geq k)).
\end{eqnarray*}
By the definition of $\{l_k(c), 1\leq k\leq N+1\}$ and $T^*$, we know that
\begin{eqnarray*}
\lim_{c\to 0}\sum_{k=2}^{k^*}\textbf{E}_{0}(v_kI(T^*\geq k))&=&0\\
\lim_{c\to \infty}\sum_{k=2}^{k^*}\textbf{E}_{0}(v_kI(T^*\geq k))&=&\sum_{k=2}^{k^*}\textbf{E}_{0}(v_k)
\end{eqnarray*}
As $\sum_{k=2}^{k^*}\textbf{E}_{0}(v_kI(T^*\geq k))$ is continuous and increasing on $c$, it follows that there is a positive number $c_{\gamma}$ such that
\begin{eqnarray}
\textbf{E}_{0}(\sum_{k=1}^{T^*(c_{\gamma})}v_k)=\sum_{k=1}^{k^*}\textbf{E}_{0}(v_kI(T^*\geq k))=\gamma.
\end{eqnarray}

It follows from (A.8) that
\begin{eqnarray}
&&\sum_{m=1}^{N}\textbf{E}_{0}([Y_m-cv_{m+1}]I(T^*\geq m+1))\nonumber\\
&=&\textbf{E}_{0}\Big(\textbf{E}_{0}(\sum_{m=1}^NB_{m, 1}(N)[Y_m-cv_{m+1}]|\mathfrak{F}_{0})\Big) \nonumber\\
&=&\textbf{E}_{0}(cv_1-l_0(c)).
\end{eqnarray}
Thus, by (A.1), (A.11), and (A.12) we have
\begin{eqnarray*}
&{}&\mathcal{J}_{M, N}(T^*(c_{\gamma}))=\frac{\textbf{E}_{0}( \sum_{m=1}^{T^*(c_{\gamma})}Y_{m-1})}{\textbf{E}_{0}(\sum_{m=1}^{T^*(c_{\gamma})}v_k)}\\
&=&\frac{ \sum_{m=1}^{N}\textbf{E}_{0}([Y_m-c_{\gamma}v_{m+1}+c_{\gamma}v_{m+1}]I(T^*(c_{\gamma})\geq m+1))}{\gamma}\\
&=&\frac{ c_{\gamma}\sum_{m=1}^{N+1}\textbf{E}_{0}(v_{m}I(T^*(c_{\gamma})\geq m))-\textbf{E}_{0}(l_0(c_{\gamma}))}{\gamma}\\
&=&\frac{c_{\gamma}\textbf{E}_{0}(\sum_{m=1}^{T^*(c_{\gamma})}v_k)}{\gamma}-\frac{\textbf{E}_{0}(l_0(c_{\gamma}))}{\gamma}=c_{\gamma}-\frac{\textbf{E}_{0}(l_0(c_{\gamma}))}{\gamma}\\
&=&c_{\gamma}(1-\frac{\textbf{E}_{0}(v_1)}{\gamma})-\frac{\textbf{E}_{0}[l_1(c_{\gamma})-Y_1]^+}{\gamma}.
\end{eqnarray*}
The last equality follows from the definition of $l_0(c)$ in (6). It proves (iii) of Theorem 1.

\textbf{Step III.} Show (i) and (ii) of Theorem 1. Let
\begin{eqnarray*}
\tilde{c}_{\gamma}= \mathcal{J}_{M, N}(T^*(c_{\gamma}))=c_{\gamma}-\frac{\textbf{E}_{0}(l_0(c_{\gamma}))}{\textbf{E}_{0}(\sum_{m=1}^{T^*(c_{\gamma})}v_m)}.
\end{eqnarray*}
If $\mathcal{J}_{M, N}(T)\geq c_{\gamma}$, then $\mathcal{J}_{M, N}(T)\geq \tilde{c}_{\gamma}= \mathcal{J}_{M, N}(T^*(c_{\gamma}))$. If $\mathcal{J}_{M, N}(T)<  c_{\gamma}$, then,
by (A.1), (A.3), and $\textbf{E}_{0}(\sum_{m=1}^{T}v_m)\geq \gamma$,  we have
\begin{eqnarray*}
[\mathcal{J}_{M, N}(T)-c_{\gamma}]\gamma &\geq & [\frac{ \textbf{E}_{0}(
\sum_{m=1}^{T}Y_{m-1})}{\textbf{E}_{0}(\sum_{m=1}^{T}v_m)}-c_{\gamma}]\textbf{E}_{0}(\sum_{m=1}^{T}v_m)\\
&=&[\textbf{E}_{0}(
\sum_{m=1}^{T}Y_{m-1})-c_{\gamma}\textbf{E}_{0}(\sum_{m=1}^{T}v_m)]\\
&\geq & [\textbf{E}_{0}( \sum_{m=1}^{T^*(c_{\gamma})}Y_{m-1})-c_{\gamma}\textbf{E}_{0}(\sum_{m=1}^{T^*(c_{\gamma})}v_m)]\\
&=&[\frac{ \textbf{E}_{0}( \sum_{m=1}^{T^*(c_{\gamma})}Y_{m-1})}{\textbf{E}_{0}(\sum_{m=1}^{T^*(c_{\gamma})}v_m)}-c_{\gamma}]\textbf{E}_{0}(\sum_{m=1}^{T^*(c_{\gamma})}v_m)\\
&=&[\mathcal{J}_{M, N}(T^*(c_{\gamma}))-c_{\gamma}]\gamma.
\end{eqnarray*}
This means that $\mathcal{J}_{M, N}(T)\geq \mathcal{J}_{M, N}(T^*(c_{\gamma}))$ for all $T\in \mathfrak{T}_N$ with $\textbf{E}_{0}(\sum_{m=1}^{T}v_m)\geq \gamma$. That is, (i) of Theorem 1  is true. The strict inequality in (ii) of Theorem 1 comes from the strict inequality in (A.3) when  $T\neq T^*(c_{\gamma})$ with $\textbf{E}_{0}(\sum_{m=1}^{T}v_m)=\textbf{E}_{0}(\sum_{m=1}^{T^*(c_{\gamma})}v_m)= \gamma$. This completes the proof of Theorem 1.

{\bf Proof of Theorem 2. } Since $Y_k=(Y_{k-1}+w_{k}(Y_{k-1}, A_{n, p_1}))\Lambda_k$ and
\begin{eqnarray*}
\Lambda_k=\frac{p_{\textbf{1}k}(X_k|X_{k-1}, ..., X_{k-j})}{p_{\textbf{0}k}(X_k|X_{k-1}, ..., X_{k-i})}
\end{eqnarray*}
for $1\leq k\leq N$, it follows that $(Y_k, X_k), 0\leq k\leq N,$ is a two-dimensional $\textbf{p}$-order Markov chain, where $p=\max\{i, j\}$.  Let $1\leq p\leq N$. By the definition of the optimal control limits, we have
\begin{eqnarray}
\\
l_k(c)&=&cv_{k+1}(Y_k, A_{n, p_2})\nonumber \\
&&+\textbf{E}_{0}\Big([l_{k+1}(c)-(Y_k+w_{k+1}(Y_k, A_{n, p_1}))\Lambda_{k+1}]^+|Y_k, A_{n, 0}\Big)\nonumber
\end{eqnarray}
for $0\leq k\leq  p-1$ and
\begin{eqnarray}
\\
l_k(c)&=&cv_{k+1}(Y_k, A_{n, p_2})\nonumber \\
&+&\textbf{E}_{0}\Big([l_{k+1}(c)-(Y_k+w_{k+1}(Y_k, A_{n, p_1}))\Lambda_{k+1}]^+|Y_k, A_{n, p}\Big)\nonumber
\end{eqnarray}
for $ p \leq k\leq N.$ Let $p=0$,  we have similarly
\begin{eqnarray}
\\
l_k(c)&=&l_k(c, Y_k)\nonumber \\
&=&cv_{k+1}(Y_k)+\textbf{E}_{0}\Big([l_{k+1}(c)-(Y_k+w_{k+1}(Y_k))\Lambda_{k+1}]^+|Y_k\Big)\nonumber
\end{eqnarray}
for $0\leq k\leq N$.

{\bf Proof of Theorem 3. }  Let $p=0$. As the observations $X_k, 0\leq k\leq N,$ are independent, it follows from the definition of $\{Y_k, 1\leq k\leq N\}$ that $\{Y_k, 1\leq k\leq N\}$ is a $1$-order Markov chain. Thus, the optimal control limits, $l_k(c), 0\leq k\leq N$, satisfy (A.15). Let $y=cv_{N+1}(y)$.  As $v_{N+1}(y)$ is non-increasing, it follows that there is a positive number $y_N(c)$ such that $y_N(c)=cv_{N+1}(y_N(c))$. Hence, $Y_N\geq l_N(c)=cv_{N+1}(Y_N)$ if and only if $Y_N\geq y_N(c)$. Therefore, we let $\widetilde{l_N}(c)=y_N(c)$. Take $k=N-1$ in (A.15) and let
\begin{eqnarray*}
y=f_0(y)&=&cv_{N}(y)\\
&+&\textbf{E}_{0}\Big([cv_{N+1}((y+w_{N}(y))\Lambda_{N})-(y+w_{N}(y))\Lambda_{N}]^+|Y_{N-1}=y\Big).
\end{eqnarray*}
Note that the two functions $(y+w_{N}(y))$ and $v_{N}(y)$ are non-decreasing and non-increasing on $y\geq 0$, respectively. Therefore, the function  $f_0(y)$ is
non-increasing on $y\geq 0$, and it follows that there is a positive number $y_{N-1}$ such that $y_{N-1}=f_0(y_{N-1})$; that is,
\begin{eqnarray*}
y_{N-1}&=&cv_{N}(y_{N-1})\\
&+&\textbf{E}_{0}\Big([cv_{N+1}-(y_{N-1}+w_{N}(y_{N-1}))\Lambda_{N}]^+|Y_{N-1}=y_{N-1}\Big).
\end{eqnarray*}
This implies that $Y_{N-1}\geq l_{N-1}(c)$ if and only if $Y_{N-1}\geq y_{N-1}$. Therefore, we let $\widetilde{l_{N-1}}(c)= y_{N-1}$. Similarly, there are  positive numbers $ y_{k},
1\leq k\leq N-2$ such that $Y_{k}\geq l_{k}(c)$ if and only if $Y_{k}\geq y_{k}$ for $1\leq k\leq N-2$, where
\begin{eqnarray*}
y_{k}=cv_{k+1}(y_k)+\textbf{E}_{0}\Big([l_{k+1}(c)-(y_{k}+w_{k+1}(y_{k}))\Lambda_{k+1}]^+|Y_{k}=y_{k}\Big)
\end{eqnarray*}
for $1\leq k\leq N-2$. Taking $\widetilde{l_k}(c)=y_{k}$ for $1\leq k\leq N$, we know that the control limit $\{ \widetilde{l_k}(c) \}$ is an equivalent control limit of the optimal sequential test $T^*_{M}(c, N)$ and it consists of a series of nonnegative non-random  numbers. This proves (ii) of Theorem 3.

Let $1\leq p\leq N.$ As $\{(Y_k, X_k), 0\leq k\leq N\}$  is a two-dimensional $p$-order Markov chain, it follows that (A.14) and (A.13) hold for  $p\leq k\leq N$ and $0\leq k\leq p-1$, respectively.
When $k=N$ in (A.14), we take  $\widetilde{l_N}(c)=y_N(c)$, where $y_N(c)=cv_{N+1}(y_N(c))$. For any fixed observation values $a_{k, p}=\{x_{k}, ..., x_{k-p+1}\}$ for $p\leq k\leq N-1$ and
$a_{k, 0}=\{x_{k},...,x_{0}\}$ for $0\leq k\leq p-1$, let
\begin{eqnarray*}
y=f_p(y)&=& cv_{k+1}(y, a_{k, p})\\
&+&\textbf{E}_{0}\Big([l_{k+1}(c)-(y+w_{k+1}(y, a_{k, p}))\Lambda_{k+1}]^+|Y_{k}=y, A_{k, p}=a_{k, p}\Big)
\end{eqnarray*}
for $p\leq k\leq N-1$ and
\begin{eqnarray*}
y=g_p(y) &=& cv_{k+1}(y, a_{k, 0} )\\
&+&\textbf{E}_{0}\Big([l_{k+1}(c)-(y+w_{k+1}(y, a_{k, 0}))\Lambda_{k+1}]^+|Y_{k}=y, A_{k, 0}=a_{k, 0}\Big)
\end{eqnarray*}
for $0\leq k\leq p-1.$ As the two functions $f_p(y)$ and $g_p(y)$ are non-increasing on $y\geq 0$, it follows that there are positive numbers $y_k=y_k(c, a_{k, p})$ for
$p\leq k\leq N-1$ and $y_k=y_k(c, a_{k, 0})$ for  $1\leq k\leq p-1$  such that $y_k=f_p(y_k)$ for $p\leq k\leq N-1$ and $y_k=g_p(y_k)$ for $1\leq k\leq p-1$. Therefore,
$Y_k\geq l_k(c)$ if and only if $Y_k\geq y_k$. Taking $\widetilde{l_k}(c)=y_k(c, X_{k},...,X_{k-p+1})$ for $p\leq k\leq N$ and  $\widetilde{l_k}(c)=y_k(c, X_{k},...,X_{0})$ for $1\leq k\leq p-1$, we have  $\widetilde{T^*_{M}}(c, N)=T^*_{M}(c, N)$. That is, $\{\widetilde{l_k}(c), \, 1\leq k\leq N+1\}$ is an equivalent control limit of the optimal sequential test $T^*_{M}(c, N)$ that does not directly depend on the statistic, $Y_k, 1\leq k\leq N.$  This completes the proof of (i) of Theorem 3.

\begin{thebibliography}{12}

\bibitem{} Akoglu, L., Tong, H. H. and Koutra, D. (2015) Graph based anormaly detection and description: a survey. \emph{Data Min Know Disc.} \textbf{29}, 626-688.

\bibitem{} Basseville, M. and Nikiforov, I. (1993) \emph{Detection of Abrupt Changes: Theory and Applications.} Prentice-Hall, Englewood Cliffs.

\bibitem{} Bersimis, S., Sgora, A. and Psarakis, S.  (2018) The application of multivariate statistical process monitoring in non-industrial processes. \emph{Qual. Technol.  Quant. Manag.}, \textbf{15}, 526-549.

\bibitem{} Bersimis,S., Psarakis, S. and Panaretos, J. (2007)  Multivariate statistical process control charts: An Overview. \emph{Qual. Reliab. Engng. Int.}, \textbf{23}, 517-543.

\bibitem{} Chakraborti, S., van der Laan, P. and Bakir, S. T. (2001) Nonparametric control charts: an overview and some results. \emph{J. Qual. Technol.}, \textbf{33}, 304-315.

\bibitem{} Chen, X. and Baron, M. (2014) Change-point alalysis of survival data with application in clinical trials.\emph{ Open J. Statistics}, \textbf{4}, 663-677.

\bibitem{} Chow, Y. S., Robbins, H.and Siegmund, D.(1971) \emph{The Theory of Optimal Stopping.} Dover Publications, INC. New York.

\bibitem{} Fris\'{e}n, M. (2003) Statistical Surveillance, Optimality and Methods. \emph{Int. Statist. Rev.}, \textbf{71}, 403-434.

\bibitem{} Fris\'{e}n, M. (2009) Optimal Sequential Surveillance for Finance, Public Health, and Other Areas (with Discussion). \emph{Sequent. Analysis}, \textbf{28}, 310-337.

\bibitem{} Han, D. and Tsung, F. G. (2004) A generalized EWMA control chart and its comparison with the optimal EWMA, CUSUM and GLR schemes.\emph{ Ann. Statist.}, \textbf{32}, 316-339.

\bibitem{}  Han, D., Tsung, F. G. and Xian, J. G. (2017) On the optimality of  Bayesian change-point detetion.  \emph{Ann. Statist.}, \textbf{45}, 1375-1402.

\bibitem{} Hosseini, S. H. and Noorossana, R. (2018) Performance evaluation of EWMA and CUSUM control charts to detect anomalies in social networks using average
and standard deviation of degree measures. \emph{Qual Reliab Engng Int.},\textbf{ 34}, 477¨C500.

\bibitem{} Lai, T. L. (1995) Sequential change-point detection in quality control and dynamical systems (with discussion). \emph{J. Roy. Statist. Soc. Ser. B.}, \textbf{57}, 613-658.

\bibitem{} Lai, T. L. (2001) Sequential analysis: some classical problems and new challenges. \emph{Statist Sinica},  \textbf{11}, 303-408.

\bibitem{} Lorden, G. (1971) Procedures for reacting to a change in distribution. \emph{Ann. Math. Statist.}, \textbf{42}, 1897-1908.

\bibitem{} Montgomery, D. C. (2009) \emph{Introduction to Statistical Quality  Control}. 6th ed. New York: John Wiley \& Sons.

\bibitem{}  Moustakides, G. V. (1986) Optimal stopping times for detecting changes in distribution. \emph{Ann. Statist.}, \textbf{14}, 1379-1387.

\bibitem{} Moustakides, G. V. (2008) Sequential change detection revisited. \emph{Ann. Statist.},\textbf{ 36}, 787-807.

\bibitem{} Page, E. S. (1954) Continuous inspection schemes. \emph{Biometrika},\textbf{ 41}  100-114.

\bibitem{} Pollak, M. (1985) Optimal detection of a change in distribution. \emph{Ann. Statist.}, \textbf{13}, 206-227.

\bibitem{} Polunchenko, A. S. and Tartakovsky, A. G. (2010). On optimality of the Shiryaev-Roberts procedure for detecting a change in distribution. \emph{ Ann. Statist.}, \textbf{38}, 3445-3457.

\bibitem{} Poor, H.V. and Hadjiliadis, O. (2009). \emph{Quickest Detection}. Cambridge University Press, Camridge, New York.

\bibitem{} Qiu, P. (2014) \emph{Introduction to Statistical Process Control}. Boca Raton, FL: Chapman \& Hall/CRC.

\bibitem{} Rigdon, S. E. and Fricker, R. D., Jr. (2015) \emph{Innovative Statistical Methods for Public Health Data} (eds by D. G. Chen and J. Wilson). ICSA Book Series in Statistics. Springer International Publishing Switzerland.

\bibitem{} Roberts, S.W. (1959) Control chart tests based on geometric moving average.  \emph{Technomietrics},  \textbf{1},  239-250.

\bibitem{} Saleh, N. A., Mahmoud, M. A., Jones-Farmer, L. A., Zwetsloot, I. and  Woodall, W. H. (2015). Another Look at the EWMA Control Chart with Estimated Parameters. \emph{J. Qual. Technol.},\textbf{ 47}, 363-382.

\bibitem{} Shewhart, W. A. (1931)  \emph{Economic Control of Quality of Manufactured Product}. New York: Van Nostrand.

\bibitem{} Shiryaev, A. N. (1963) On optimum methods in quickest detection problems. \emph{Theory Probab. Appli.}, \textbf{13}, 22-46

\bibitem{} Shiryaev, A. N. (1978) \emph{Optimal Stopping Rules}. Springer-Verlag, New York.

\bibitem{} Siegmund, D. (1985) \emph{Sequential Analysis: Tests and Confidence Intervals}. Springer Series in Statistics.
Springer-Verlag, New York

\bibitem{} Siegmund, D. (2013) Change-points: from sequential detection to biology and back.\emph{ Sequ.  Analysis}, \textbf{32}, 2-14

\bibitem{} Stoumbos, Z. G., Reynolds, M. R., Ryan, T. P., and Woodall, W. H. (2000) The state of statistical process control as we proceed into the 21st century. \emph{J. Amer. Statist. Assoc.}, \textbf{95}, 992-998.

\bibitem{} Tartakovsky, A. G., Nikiforov, I. and Basseville, M, (2015).  \emph{Sequential Analysis:  Hypothesis Testing and  Changepoint Detection}. CRC Press, Taylor $\&$  Francis Group.

\bibitem{} Woodall, W. H. (2006) The use of control charts in health-care and public-health survellance. \emph{J. Qual. Technol.},\textbf{ 38}, 89-118.

\bibitem{} Woodall, W. H. and Montgomery, D. C. (2014) Some current directions in the theory and application of statistical process monitoring. \emph{J. Qual. Technol.},\textbf{ 46}, 78-94.

\bibitem{} Woodall, W. H., Zhao, M. J., Paynabar, K., Sparks, R. and Wilson, J. D. (2017)  An overview and perspective on social network monitoring. \emph{IIE Trans.}, \textbf{49}, 354-365.


\end{thebibliography}
\end{document}